\documentclass[12pt]{article}
\usepackage{amsmath,amsfonts,amssymb,amsthm, amstext}

\newtheorem{proposition}{Proposition}[section]
\newtheorem{theorem}[proposition]{Theorem}
\newtheorem{corollary}[proposition]{Corollary}
\newtheorem{lemma}[proposition]{Lemma}
\newtheorem{remark}[proposition]{Remark}

\newcommand{\nc}{\newcommand}
\nc{\R}{{\mathbb R}}
\nc{\N}{{\mathbb N}}
\nc{\Z}{{\mathbb Z}}

\nc{\BP}{\mathbb{P}}
\nc{\BQ}{\mathbb{Q}}
\nc{\BE}{\mathbb{E}}

\begin{document}

\author{{\sc J\"urgen Kampf\footnotemark, G\"unter Last\footnotemark}
} 
\title{On the convex hull of symmetric stable processes}
\date{\today}
\maketitle
\setcounter{footnote}{1} \footnotetext{Postal address: Institut f\"ur Stochastik, Karlsruhe Institute of Technology,
76128 Karlsruhe, Germany. 
Email address: juergen.kampf@kit.edu} 
\stepcounter{footnote} \footnotetext{Postal address: Institut f\"ur Stochastik, Karlsruhe Institute of Technology,
76128 Karlsruhe, Germany. 
Email address: guenter.last@kit.edu} 
\begin{abstract} 
\noindent 
Let $\alpha\in (1,2]$ and $X$ be an $\R^d$-valued $\alpha$-stable
process with independent and symmetric components starting
in $0$. We consider the closure $S_t$ of the path described by $X$
on the interval $[0,t]$ and its convex hull $Z_t$.
The first result of this paper provides a formula for certain mean
mixed volumes of $Z_t$ and in particular for the
expected first intrinsic volume of $Z_t$.
The second result deals with the asymptotics
of the expected volume of the stable sausage
$Z_t+B$ (where $B$ is an arbitrary convex body with interior points) as $t\to 0$. 
\end{abstract}


\vspace{0.2cm}

\section{Introduction and main results}\label{intro}
\setcounter{equation}{0}

For fixed $\alpha\in(1,2]$ and fixed integer $d\ge 1$ 
we consider an $\R^d$-valued stochastic process 
$X\equiv (X(t))_{t\ge 0}=(X_1(t),\ldots,X_d(t))_{t\ge 0}$,
defined on the probability space
$(\Omega,\mathcal{A},\BP)$, such that the components
$X_j:=(X_j(t))_{t\ge 0}$, $j\in\{1,\ldots,d\}$, 
are independent $\alpha$-stable symmetric L\'evy processes
with scale parameter $1$ starting in $0$.
The characteristic function of $X_j(t)$ is given by
\begin{align}\label{chf}
\BE\exp[isX_j(t)]=\exp[-t|s|^\alpha],\quad s\in\R,t\ge 0,
\end{align}
cf.\ \cite[Section 1.3]{SaTa94}.
This implies that $X$ is {\em self-similar} in the sense that
$(X(st))_{s\ge 0}\overset{d}{=}t^{1/\alpha}X$ for any $t>0$,
see \cite[Example 7.1.3]{SaTa94} and \cite[Chapter 15]{Kallenberg}.
We assume that $X$ is right-continuous with
left-hand limits (rcll). 
For $t\ge 0$, let $S_t$ be the closure of the path 
$S^0_t:=\{X(s):0\le s\le t\}$
and let $Z_t$ denote the convex hull of $S_t$.
These are random closed sets. We abbreviate $Z:=Z_1$.
By self-similarity
\begin{align}\label{selfsim}
Z_t\overset{d}{=}t^{1/\alpha}Z,\quad t>0.
\end{align}

If $\alpha=2$ then $X$ is a standard Brownian motion.
A classical result of \cite{Takacs80} for planar Brownian
motion says that
\begin{align}\label{1.1}
\BE V_1(Z)=\sqrt{2\pi},
\end{align}
where $V_1(K)$ denotes half the circumference of a 
convex set $K\subset\R^2$.
Our first aim in this paper is to formulate and to
prove such a result for arbitrary $\alpha\in (1,2]$ and
arbitrary dimension $d$. 
In fact we also consider more general geometric functionals.
A {\em convex body} (in $\R^d$) is a non-empty compact and convex
subset of $\R^d$. 
We let $V(K_1,\ldots,K_d)$, denote the mixed volumes of convex bodies 
$K_1,\ldots,K_d\subset \R^d$ \cite[Section 5.1]{Schneider1993}.
These functionals are symmetric in $K_1,\ldots,K_d$
and we have for any convex bodies $K,B\subset\R^d$
\begin{align}\label{genS}
V_d(K+tB)=\sum^d_{j=0}\binom{d}{j}t^{d-j}V(K[j],B[d-j]),\quad t\ge 0,
\end{align}
where $V_d$ is Lebesgue measure, $tB:=\{tx:x\in B\}$,
$B+C:=\{x+y:x\in B,y\in C\}$ is the {\em Minkowski sum}
of two sets $B,C\subset\R^d$, and 
$V(K[j],B[d-j])$ is the mixed volume of $K_1,\ldots,K_d$ in case
$K_1=\ldots=K_j=K$ and $K_{j+1}=\ldots=K_d=B$. The $j$th
{\it intrinsic volume} $V_j(K)$ of a convex body $K$ is given by
\begin{align}\label{v1}
V_j(K)=\frac{\binom{d}{j}}{\kappa_{d-j}}V(K[j],B^d[d-j]),\quad j=0,\ldots,d,
\end{align}
where $B^d$ is the Euclidean unit ball in $\R^d$, and
$\kappa_j$ is the $j$-dimensional volume of $B^j$.  
In particular, $V_d(K)$ is the volume 
of $K$, $V_{d-1}(K)$ is half the surface area, $V_{d-2}(K)$ is proportional to 
the integral mean curvature, $V_1(K)$ 
is proportional to the mean width of $K$, and $V_0(K)=1$.
(If $d=2$ then $V_1(K)$ has been discussed at \eqref{1.1}.) 
A geometric interpretation of $V(K_1,\ldots,K_d)$ in the case
$K_1=\ldots=K_{d-1}=B$ is provided by \eqref{genS}:
\begin{align}\label{1.2}
V(B,\ldots,B,K)=\lim_{r\to 0}r^{-1}(V_d(B+rK)-V_d(B)).
\end{align}

For any $p\ge 1$ define $B_p:=\{u\in\R^d:\|u\|_p\le 1\}$ as the unit ball
with respect to the $L_p$-norm $\|(u_1,\ldots,u_d)\|_p:=(|u_1|^p+\ldots+|u_d|^p)^{1/p}$.
Finally we introduce the constant
\begin{align}\label{ca}
c_\alpha:=\frac{\alpha}{2}\BE|X_1(1)|,
\end{align}
Since $\alpha>1$, this constant is finite, see \cite[Property 1.2.16]{SaTa94}. 
A direct calculation shows that 
\begin{align}\label{c2}
c_2=\sqrt{\frac{2}{\pi}}.
\end{align}
In the case of $\alpha<2$ we are not aware of an explicit
expression for $c_\alpha$.

\begin{theorem}\label{t1} Let $K_1,\ldots,K_{d-1}\subset \R^d$ be 
convex bodies. Then
\begin{align}\label{main1}
\BE V(K_1,\ldots,K_{d-1},Z)=c_\alpha V(K_1,\ldots,K_{d-1},B_{\alpha'}),
\end{align}
where $1/\alpha+1/\alpha'=1$.
\end{theorem}

\begin{remark}\rm By the scaling relation \eqref{selfsim} and
the homogeneity property of mixed volumes \cite[(5.1.24)]{Schneider1993}
the identity \eqref{main1} can be generalized to
\begin{align}\label{main1a}
\BE V(K_1,\ldots,K_{d-1},Z_t)=
c_\alpha t^{1/\alpha}V(K_1,\ldots,K_{d-1},B_{\alpha'}).
\end{align}
A similar remark applies to all results of this paper.\end{remark}

The proof of Theorem \ref{t1} relies on the fact that
\begin{align}\label{meanbody}
\BE h(Z,u)=h(B_{\alpha'},u),\quad u\in S^{d-1},
\end{align}
where $S^{d-1}$ denotes the unit sphere,
\begin{align*}
h(K,u):=\max\{\langle x,u\rangle:x\in K\},\quad u\in S^{d-1},
\end{align*} 
is the {\em support function} of a convex body $K$, and $\langle\cdot,\cdot\rangle$
denotes the Euclidean scalar product on $\R^d$.
This means that $B_{\alpha'}$ is the {\em mean body} of $Z$
\cite[p.\ 146]{SW}, or the {\em selection expectation} of $Z$ 
\cite[Theorem 2.1.22]{Molchanov05}.

The next corollary provides a direct generalization of
\eqref{1.1}.

\begin{corollary}\label{rbrown} Assume that $X$ is a standard-Brownian motion in $\R^{d}$. Then
\begin{align}\label{V1S}
\BE V_1(Z)=
\frac{d\sqrt{2} \Gamma(\frac{d-1}{2}+1)}{\Gamma(\frac{d}{2}+1)}
\end{align}
\end{corollary}

In the case of Brownian motion it is possible to
calculate the expectation of the second intrinsic volume
$V_2(Z)$ of $Z$.

\begin{proposition}\label{p1} 
Assume that $X$ is a standard-Brownian motion in $\R^{d}$. Then
$$
\BE V_2(Z)=(d-1)\frac{\pi}{2}.
$$
\end{proposition}

Our second theorem deals with the asymptotic behaviour of the expected volume of the
{\em stable sausage} $S_t+B$ as $t\to 0$, where
$B$ is a convex body.
Our result complements classical results 
on the asymptotic behaviour 
of $\BE V_d(S_t+B)$ as $t\to \infty$, cf.\  \cite{Spitzer} for the case of
Brownian motion and  \cite{Getoor65} for the case of more
general stable processes.

\begin{theorem}\label{t2} 
Let $B$ be a convex body with non-empty interior
and let $\alpha'$ be as in Theorem \ref{t1}. Then
$$
\lim_{t\to 0} t^{-1/\alpha}(\BE V_d(S_t+B)-V_d(B))
=d c_\alpha V(B,\ldots,B,B_{\alpha'}).
$$
\end{theorem}

In the case $\alpha=2$ the random set $S_t+B$ is
known as {\em Wiener sausage}. Even then Theorem \ref{t2} seems to be new:

\begin{corollary}\label{co2} Assume that $X$ is a
Brownian motion and let $B$ be a convex body with non-empty interior.
Then
$$
\lim_{t\to 0} t^{-1/2}(\BE V_d(S_t+B)-V_d(B))
=\frac{d\sqrt{2}}{\sqrt{\pi}} V(B,\ldots,B,B^d).
$$
In the case $B=B^d$ the limit equals 
$2\sqrt{2}\pi^{(d-1)/2}/\Gamma(d/2)$.
\end{corollary}

\begin{remark}\rm In the special case $d=3$ and $\alpha=2$ we have (see \cite{Spitzer})
\begin{equation}\mathbb{E}\,V_3(S_t+rB^3)=\frac{4}{3}\pi r^3 + 4\sqrt{2\pi} r^2\sqrt{t} + 2\pi rt\label{3vol} \end{equation}
for any $r,t\ge 0$. The term constant in $t$ clearly allows a geometric interpretation as $V_3(rB^3)$. Now we are able to give a geometric interpretation of the coefficient of $\sqrt{t}$ as well. 

From \eqref{3vol} we get
$$\lim_{t\to 0} t^{-1/2}(\BE V_3(S_t+rB^3)-V_3(rB^3))=\lim_{t\to 0}t^{-1/2}(4\sqrt{2\pi} r^2\sqrt{t} + 2\pi rt)=4\sqrt{2\pi} r^2.$$
On the other hand from the proof of Theorem \ref{t2} one can see
$$\lim_{t\to 0} t^{-1/2}(\BE V_3(S_t+rB^3)-V_3(rB^3))=3\BE V(rB^3,rB^3,Z).$$
By \eqref{v1} and the homogenity property of mixed volumes (see e.g. \cite[(5.1.24)]{Schneider1993}) we have
$$3 V(rB^3,rB^3,Z)=r^2\kappa_2V_1(Z).$$
Altogether this is
$$4\sqrt{2\pi} r^2=r^2\kappa_2\BE V_1(Z).$$
\end{remark}

\section{Proofs}\label{sec2}
\setcounter{equation}{0}

We need the following measurability property
of the closure $S_t$ of $\{X(s):0\le s\le t\}$ and its convex hull $Z_t$,
refering to \cite{Molchanov05,SW} for the notion
of a {\em random closed set}.

\begin{lemma}\label{l1} For any $t\ge 0$, $S_t$ and $Z_t$
are random closed sets.
\end{lemma}
{\sc Proof:} To prove the first assertion it is enough
to show that $\{S_t\cap G=\emptyset\}$ is measurable
for any open $G\subset\R^d$, see \cite[Lemma 2.1.1]{SW}.
But since $X$ is rccl it is clear that
$S_t\cap G=\emptyset$ iff $X(u)\notin G$
for all rational numbers $u\le t$.
The second assertion is implied by \cite[Theorems 12.3.5,12.3.2]{SW}.\qed

\vspace*{0.3cm}
The previous lemma implies, for instance, that
$V(K_1,\ldots,K_{d-1},Z_t)$ and $V_d(S_t+B)$ are random variables,
see e.g.\ \cite[p.\ 275]{Schneider1993} and \cite[Theorem 12.3.5 and Theorem 12.3.6]{SW}.

\vspace*{0.3cm}
{\sc Proof of Theorem \ref{t1}:}
By \cite[(5.1.18)]{Schneider1993} we have that
\begin{align}\label{5118} 
V(K_1,\ldots,K_{d-1},K)=\frac{1}{d}\int_{S^{d-1}} h(K,u)\, 
S(K_1,\ldots,K_{d-1},du) 
\end{align}
holds for every convex body $K\subset \R^d$, where 
$S(K_1,\ldots,K_{d-1},\cdot)$ is the {\em mixed area measure} 
of $K_1,\ldots,K_{d-1}$, 
see \cite[Section 4.2]{Schneider1993}.  
From \eqref{5118} and Fubini's theorem we obtain that
\begin{align}\label{2.5} 
\BE V(K_1,\ldots,K_{d-1},Z)=\frac{1}{d}\int_{S^{d-1}} \BE h(Z,u)\, 
S(K_1,\ldots,K_{d-1},du).
\end{align}
For any $u\in S^{d-1}$ we have
\begin{align*}
\BE h(Z,u) &=\BE\max\{\langle x,u\rangle: x\in Z_1\}\\ 
&=\BE\sup\{\langle x,u\rangle: x\in S^0_1\}\\ 
&=\BE\sup\{\langle X(s),u\rangle: s\in [0,1]\}.
\end{align*}
It follows directly from \eqref{chf} that the process $\langle X,u\rangle$
has the same distribution as $\|u\|_\alpha X_1$.
By \cite[Theorem 4a]{Bingham73}, $\sup\{X_1(s):s\in [0,1]\}$
has a finite expectation.
Differentiating equation (7b) in \cite{Bingham73} (Spitzer's identity in continuous
time), one can easily show that
$$
\BE \sup\{X_1(s):s\in [0,1]\}=\alpha\BE X_1(1)^+,
$$
where $a^+:=\max\{0,a\}$ denotes the positive part of a real
number $a$. Since $X_1(1)$ has a symmetric distribution
and $\BP(X_1(1)=0)=0$ (stable distributions have a density) 
we have $\BE|X_1(1)|=2\BE X_1(1)^+$
and it develops that $\BE h(Z,u)=c_\alpha\|u\|_\alpha$, with $c_\alpha$ given 
by \eqref{ca}. Inserting this result into \eqref{2.5} gives
\begin{align}\label{2.9}
\BE V(K_1,\ldots,K_{d-1},Z)=\frac{c_\alpha}{d}\int_{S^{d-1}} 
\|u\|_\alpha\, S(K_1,\ldots,K_{d-1},du).
\end{align}
By \cite[Remark 1.7.8]{Schneider1993}, $\|u\|_\alpha$ is the support
function of the {\em polar body}
$$
B^*_\alpha:=\{x\in\R^d:\text{$\langle x,u\rangle\le 1$ for all $u\in B_\alpha$}\}
$$
of $B_\alpha$. Using the H\"older inequality, it is straightforward
to check that $B^*_\alpha=B_{\alpha'}$, where $1/\alpha+1/\alpha'=1$.
Using this fact as well as \eqref{5118} in \eqref{2.9}, we obtain
the assertion \eqref{main1}.\qed

\vspace{0.3cm}
{\sc Proof of Corollary \ref{rbrown}:} 
Since $\alpha=2$ we have $\alpha'=2$ and $B_{\alpha'}=B^d$.
By Theorem \ref{t1} and \eqref{v1},
\begin{align}\label{2.91}
\BE V_1(Z)=\BE \frac{d}{\kappa_{d-1}} V(B^d,\dots,B^d,Z)=\frac{c_2d}{\kappa_{d-1}}V(B^d,\ldots,B^d, B^d)
=\frac{c_2d\kappa_d}{\kappa_{d-1}},
\end{align}
where we have used that $V(B^d,\ldots, B^d)=V_d(B^d)$.
Using \eqref{c2} and the well-known formula 
$\kappa_d=\pi^{d/2}/\Gamma(d/2+1)$ in \eqref{2.91}, we obtain
the result.\qed

\vspace{0.3cm}
{\sc Proof of Proposition \ref{p1}:} 
By Kubota's formula (see e.g.\ \cite[(5.3.27)]{Schneider1993}) we have
$$
V_2(Z)=\frac{d(d-1)\kappa_d}{2\kappa_2\kappa_{d-2}}
\int_{G_{2}}V_2(Z|L)\, \nu_2(dL),
$$
where $G_{2}$ denotes the set of all $2$-dimensional linear subspaces of 
$\R^{d}$, $\nu_2$ is the Haar measure on $G_2$ with $\nu_2(G_{2})=1$ 
and $Z|L$ denotes the image of $Z$ under the orthogonal projection 
onto the linear subspace $L$. By Fubini's theorem,
\begin{align*}
\BE V_2(Z)
=\frac{d(d-1)\kappa_d}{2\kappa_2\kappa_{d-2}}
\int_{G_{2}}\BE V_2(Z|L)\, \nu_2(dL).
\end{align*}
The spherical symmetry of Brownian motion implies that $\BE V_2(Z|L)$
does not depend on $L$. Assume that
$L=\{(x_1,x_2,0,\dots,0):x_1,x_2\in\mathbb{R}\}$.
Now it is clear from the definition of the $d$-dimensional Brownian motion
that the random closed set $Z|L$ is the 
convex hull of a Brownian path in $L$. 
By Remark (a) in \cite[p.\ 149]{CHM89} (see also \cite{MaCoRa10})
we have $\BE V_2(Z|L)=\pi/2$. 
Therefore,
\begin{align*}
\BE V_2(Z)
=\frac{d(d-1)\kappa_d}{2\kappa_2\kappa_{d-2}}\frac{\pi}{2}
\end{align*}
and the result follows by a straightforward calculation.\qed


\vspace{0.3cm} 
{\sc Proof of Theorem \ref{t2}:}
By self-similarity and the dominated convergence theorem, 
on whose conditions we will comment below, we have
\begin{align} 
\lim_{t\to 0} t^{-1/\alpha}(\BE V_d(S_t+B)-V_d(B))
&=\lim_{t\to 0} t^{-1/\alpha}(\BE V_d(t^{1/\alpha}S_1+B)-V_d(B))\notag\\
&=\lim_{t\to 0} t^{-1}(\BE V_d(tS_1+B)-V_d(B))\notag\\
&=\BE\lim_{t\to 0} t^{-1}(V_d(tS_1+B)-V_d(B)).\label{f26}
\end{align}
In order to justify the application of the dominated convergence theorem, put
$$
Y_j=\sup\{X_j(s):s\in[0,1]\},\qquad
\tilde Y_j=\inf\{X_j(s):s\in[0,1]\},\quad j=1,\dots, d.
$$
As noted in the proof of Theorem \ref{t1}, $Y_j$ has a finite expectation.
Since $-\tilde Y_j$ has the same distribution as $Y_j$,
$\tilde Y_j$ has a finite expectation as well. 
From \eqref{genS} we obtain for all $t\in (0,1]$ that
\begin{align*}
t^{-1}V_d(tS+B)-V_d(B)&\le t^{-1}(V_d(tZ+B)-V_d(B))\\
&=\sum_{j=0}^{d-1} \binom{d}{j}t^{d-j-1}V(B[j],Z[d-j])\\
&\le \sum_{j=0}^d \binom{d}{j}V(B[j],Z[d-j])\\
&= V_d(Z+B). 
\end{align*}
Furthermore,
$$
Z+B\subset\mathop{\times}_{j=1}^d [\tilde Y_j-h_B(-e_j), Y_j+h_B(e_j)],
$$
where $e_j$ denotes the $j$th unit vector.
It follows that
$$
t^{-1}(V_d(tS_1+B)-V_d(B))
\le \prod_{j=1}^d \left(Y_j+h_B(e_j)-\tilde Y_j+h_B(-e_j)\right),\quad t\in (0,1].
$$
This is a product of independent random variables with finite 
expected values and hence has finite expected value.
 
By \cite[Corollary 3.2 (2)]{KiRa06} we have
$$
\lim_{t\to 0} t^{-1}(V_d(tS_1+B)-V_d(B))= 
\int_{S^{d-1}}h(Z,u)\,S(B,\ldots,B,du), 
$$
and using Theorem \ref{t1} we conclude from \eqref{f26}
\begin{align}\notag 
\lim_{t\to 0} t^{-1/\alpha}\BE (V_d(S_t+B)-V_d(B))
&=\BE\int_{S^{d-1}}h(Z,u)\,S(B,\ldots,B,du)\\\notag
&=d\BE V(B,\ldots,B,Z)\\
&=d c_\alpha V(B,\ldots,B,B_{\alpha'}).
\tag*{\qed}
\end{align}

\end{document}